\documentclass[12pt, reqno]{amsart}
\usepackage{amsmath, amsthm, amscd, amsfonts, amssymb, graphicx, color}
\usepackage[bookmarksnumbered, colorlinks, plainpages]{hyperref}

 \newtheorem{theorem}{Theorem}[section]

\newtheorem{proposition}[theorem]{Proposition}

\theoremstyle{definition}
\newtheorem{definition}[theorem]{Definition}
\newtheorem{example}[theorem]{Example}

\theoremstyle{remark}

\numberwithin{equation}{section}

\begin{document}
\setcounter{page}{1}


\title[A new class of ideal Connes amenability ]{A new class of ideal Connes amenability }

\author[A. Minapoor and A. Rejali]{$^1$Ahmad Minapoor, $^2$ Ali Rejali $^{*}$}

 \address{$^1$ Department of Mathematics, Ayatollah Borujerdi University, Borujerd, Iran}
\address{$^2$ Department of Pure Mathematics, Faculty of Mathematics and Statistics, University of Isfahan, Isfahan 81746-73441, Iran }
\email{$^1$\textcolor[rgb]{0.00,0.00,0.84}{shp\_np@yahoo.com }}
\email{$^2$ rejali@sci.ui.ac.ir }

\date{Received: xxxxxx; Revised: yyyyyy; Accepted: zzzzzz.
\newline \indent $^{*}$ Corresponding author}

\subjclass[2000]{46H20, 46H25.}

 \keywords{dual Banach algebra; amenability; Ideal Connes-amenable;.}


\maketitle

 \begin{abstract}
In this paper, we introduce a  new notion of amenability, $\sigma-$Connes ideal, say, for a large class of  dual Banach algebras. We extend the concept of ideal Connes amenability and study their properties.  Let $\sigma$ be a
$weak^{*}$-continuous endomorphism on a dual Banach algebra $\mathcal{A}$ with dense range. Then
the concept of ideal Connes-amenability and $\sigma-$ ideally Connes amenability are the same. We gave some general results and hereditary properties with some examples for this new notion
of amenability.

 \end{abstract}

\section{Introduction}
\noindent The study of the cohomologies of Banach algebras is an active and important area of research in Banach algebras. Several aspects of these cohomologies have been researched into by several researchers leading to the introduction of various notions of amenability. For example,
Johnson in \cite{joh} showed that the locally compact group $G$ is amenable as a group if and only if the group algebra
$L^{1}(G)$ is amenable as a Banach algebra. The amenability of $L^{1}(G)$ as a Banach algebra is equivalent to saying
that $L^{1}(G)$ has a vanishing first order Hochschild cohomology with coefficients in dual Banach $L^{1}(G)-$ bimodule. For a general Banach
algebra $\mathcal{A},$ Johnson in \cite{joh} showed that the amenability of $\mathcal{A}$ is equivalent to $\mathcal{A}$ having a virtual diagonal, which is also equivalent to $\mathcal{A}$ having a bounded approximate diagonal. For details on this, other important results and
important references, see \cite{joh,m3,m4,m,m1,MS} and the references contained therein.

 \noindent The dual space structure of the von Neumann algebras and the $w^*$-topology on dual Banach algebra $\mathcal{A}$ is very important
and have to be put into consideration when studying cohomologies of dual Banach algebra $\mathcal{A}.$ Thus, whenever an algebra carries a canonical $w^*$-topology, it is natural to consider derivations which are $w^*$-continuous. Thus, the notion of amenability that fits into
dual Banach algebras and for von Neumann algebras was introduced by Johnson et. al. in \cite{jkr} for von Neumann algebras,
then by Helemeskii \cite{hel} and Runde in \cite{run1} for larger class of dual Banach algebras. This notion of amenability is called
Connes amenability. A dual Banach algebra $\mathcal{A}$ is Connes-amenable if, for every normal, dual Banach $\mathcal{A}.$-bimodule $X$, every $weak^{*}$-continuous derivation $D\in \mathcal
Z^1(\mathcal{A},X)$ is inner. Clearly, the notion of Connes amenability fully exploits the important dual space structure and the $w^*$-topology on dual Banach algebra $\mathcal{A}.$ For further details on this, see \cite{hel,jkr,run1,run3}.

 \noindent Furthermore, as a continuation in the study of cohomologies of Banach algebra, Gordji and Yazdanpanah in \cite{ey} introduced two notions of amenability for a Banach algebra $A.$ These new notions are
the concepts of $I$-weak amenability and ideal amenability for Banach algebras, where $I$ is a closed ideal in $A.$ They gave some examples
to show that the notion of ideal amenability is different and weaker than the notion of weak amenability and established some general
results on the new notions of amenability. Moreover, they obtained some connections between ideal amenability and weak amenability
and concluded by posing some open questions. One of these open questions was answers partially by Mewomo in \cite{m2}.
An alternative notion which is related to ideal amenability was investigated in \cite{tbe}

 \noindent We noted and pointed out that the notion of ideal amenability introduced in \cite{ey} does not fit into dual Banach algebra due to the
dual space structure and the $w^*$-topology on dual Banach algebra. To circumvent this weakness and as a further generalization of the notion of ideal amenability, the first author, Bodaghi and Ebrahimi Bagha in \cite{mbe} introduced the notion of ideal Connes amenability for dual Banach algebras, which is weaker than the notion of Connes-amenability. For example they show that von Neumann algebras are always ideal Connes amenable.

 \noindent Motivated by the above results and the ongoing research interest in this direction, we introduce $\sigma-$ideal Connes amenability for dual Banach algebra $\mathcal{A}$,
where $\sigma$ is a $weak^{*}$-continuous endomorphism of $\mathcal{A}$, that extends the notion of ideal Connes-amenability for a large class of dual Banach algebras.

 \noindent The organization of the paper is as follows: Section \ref{Se2} contains some basic and useful definitions and the two new notions
of amenability that are introduced in this work. These are needed in subsequent sections of this study. In Section \ref{Se3},
we give some general theory and hereditary properties and establish the condition under which
the ideal Connes-amenability and $\sigma-$ ideally Connes amenability are equivalent for dual Banach algebra $\mathcal{A}.$
Some useful and important examples to illustrate our main results are given in Section \ref{Se4}. We then conclude in Section \ref{Se5}.

\section{Preliminaries and definitions}\label{Se2}
\noindent In this section, we recall some basic concepts and standard definitions that are needed in subsequent sections. We
also introduce our new notions of amenability that are studied in this work. Throughout this work, $\mathcal{A}$ will denote a
dual Banach algebra. For more details and other basic definitions, standard notions and relevant results, see \cite{mbe} and the references contained therein.

 \begin{definition}
 
 A Banach algebra is dual if there is a closed submodule $(\mathcal{A_{*})}$ of $\mathcal{A^{*}}$ such that $(\mathcal{A_{*})^{*}=\mathcal{A}}$.
Let $\mathcal{A}$ be a dual Banach algebra, $(\mathcal{I_{*})^{*}=\mathcal{I}}$ be a $w^{*}$-closed two-sided ideal in $\mathcal{A}$ and $\sigma$ be a $w^{*}$-continuous endomorphism of $\mathcal{A}$, i.e. a $w^{*}$-continuous algebra homomorphism from $\mathcal{A}$ into $\mathcal{A}$. We say that $\mathcal{I}$ has the $\sigma-$dual trace extension property if for each $\lambda \in \mathcal{I}$ with $\sigma(a) \cdot \lambda = \lambda \cdot \sigma(a) $ $(a \in \mathcal{A}),$ there exist a $\tau \in \mathcal{A}$ such that $\tau \vert_{\mathcal{I_{*}}}= \lambda $ and $a\cdot \tau- \tau \cdot a = 0$ for all $a\in \mathcal{A}$.
\end{definition}

Let $\mathcal{A}$ be a dual Banach algebra. A $w^{*}$-continuous $\sigma-$derivation from $\mathcal{A}$ into a dual Banach $\mathcal{A}$-bimodule $X$, is a $w^{*}$-continuous map
$D:A\longrightarrow X$ satisfying $$D(ab)=D(a)\cdot \sigma(b)+\sigma(a)\cdot D(b) \qquad (a,b \in A).$$
For each $x\in X$, the mapping $\delta_{x}^{\sigma}:A\longrightarrow X $, defined by $\delta_{x}^{\sigma}(a)= \sigma(a) \cdot x-x\cdot \sigma(a)$, for all $a\in \mathcal{A}$, is
a $\sigma-$derivation called an inner $\sigma-$derivation.

 \begin{definition}
A dual Banach algebra $\mathcal{A}$ is said to be $\sigma-$Connes amenable, if every $w^{*}$-continuous $\sigma-$derivation
from $\mathcal{A}$ into every normal dual $\mathcal{A}$-bimodule is an inner $\sigma-$derivation.

Let $\mathcal{A}$ be a dual Banach algebra, $\mathcal{I}$ be a $w^{*}$-closed
two-sided ideal in $\mathcal{A}$ and $\sigma$ a $w^{*}$-continuous endomorphism on $\mathcal{A}.$
We say that
\begin{itemize}
\item $\mathcal{A}$ is $\sigma-\mathcal{I}-$Connes
amenable, if every $w^{*}$-continuous $\sigma-$derivation from $\mathcal{A}$ into $\mathcal{I}$ is an inner $\sigma-$derivation.
\item $\mathcal{A}$ is $\sigma$-ideally Connes amenable, if
it is $\sigma-\mathcal{I}-$Connes amenable, for every $w^{*}$-closed two-sided ideal $\mathcal{I}$ in $\mathcal{A}$.
\end{itemize}

 \end{definition}

 \noindent Let $\sigma$ be the identity map on $\mathcal{A}$. Then, it is obvious that every
ideal Connes amenable dual Banach algebra is $\sigma$-ideally Connes
amenable. We denote $\mathcal Z^1_{\sigma,w^*}(\mathcal{A},X^*)$ for
the $w^*$-continuous $\sigma-$derivations from $A$ into $X^*$ and
denote the space of all inner $\sigma-$derivations from
$\mathcal{A}$ into $X^*$ by $\mathcal
N^1_{\sigma}(\mathcal{A},X^*)$. In addition, $\mathcal
H^1_{\sigma,w^*}(\mathcal{A},X^*)=\mathcal
Z^1_{\sigma,w^*}(\mathcal{A},X^*)/\mathcal
N^1_{\sigma}(\mathcal{A},X^*)$.

 \section{Main results}\label{Se3}
\noindent In this section, we present and prove our main results.
\begin{theorem}\label{00}
Let $\mathcal{A}$ be a dual Banach algebra and $\sigma$ be a $w^{*}-$continuous endomorphism of $\mathcal{A}$ with a dense range. Then
the ideal Connes amenability and $\sigma$-ideally Connes amenability are equivalent.
\end{theorem}

 \begin{proof}
\noindent Suppose that $\mathcal{A}$ is $\sigma$-ideally Connes amenable. Let $\mathcal{I}$ be a $w^{*}$-closed ideal of $\mathcal{A}$ and $D: \mathcal{A} \longrightarrow \mathcal{I}$ be a $w^{*}-$continuous derivation. It is easy to see that $Do\sigma: \mathcal{A} \longrightarrow \mathcal{I}$ is a $w^{*}-$continuous $\sigma-$derivation. By the assumption, there exists $x\in \mathcal{I}$ such that $D(\sigma(a))=\sigma(a)x-x\sigma(a)$, $a\in \mathcal{A}$. For an arbitrary element $b\in \mathcal{A}$, there exists $(a_{i})_{i}\subset \mathcal{A}$, such that $b=lim_{i}\sigma(a_{i})$. Hence $D(b)=lim_{i}D(\sigma(a_{i}))= lim_{i}\sigma(a_{i})x-x\sigma(a_{i})= bx-xb$. Thus $D=ad_{x}$, as required. The converse is immediate.
\end{proof}

\begin{theorem} \label{them0}
\noindent Let $\mathcal{A}$ be $\sigma$-ideally Connes amenable dual Banach algebra, where $\sigma$ is a $w^{*}-$continuous endomorphism of $\mathcal{A}$, and $\mathcal{I}$ be $w^{*}$-closed two-sided ideal in $\mathcal{A}$ with the $\sigma$-dual trace extension property. Then $\dfrac{\mathcal{A}}{\mathcal{I}}$ is $\sigma$-ideally Connes amenable dual Banach algebra.
\end{theorem}

 \begin{proof}
\noindent Let $\dfrac{\mathcal{J}}{\mathcal{I}}$ be a $w^{*}$-closed two-sided ideal in $\dfrac{\mathcal{A}}{\mathcal{I}}$. Then $\mathcal{J}$ is a $w^{*}$-closed two-sided ideal in $\mathcal{A}$, also $^{\bot}\mathcal{I}$ is a predual of $\dfrac{\mathcal{A}}{\mathcal{I}}$. We know that $ ^{\bot}\mathcal{I}$ is a closed $\mathcal{A}$-submodule of $\mathcal{J_{*}}$. Let $\pi_{*}: \mathcal{J_{*}} \longrightarrow $ $^{\bot}\mathcal{I}$ be the natural projection $\mathcal{A}$-bimodule morphism and $q: \mathcal{A} \longrightarrow \dfrac{\mathcal{A}}{\mathcal{I}}$ be the natural quotient map and $(\pi_{*})^{*}$ be the adjoint of $\pi_{*}$. Let $D: \dfrac{\mathcal{A}}{\mathcal{I}}\longrightarrow \dfrac{\mathcal{J}}{\mathcal{I}}$ be a $w^{*}$-continuous $\sigma$-derivation. Then $d=(\pi_{*})^{*}\circ D\circ q: \mathcal{A}\longrightarrow \mathcal{J}$ is a $w^{*}$-continuous $\sigma$-derivation. Let $j_{*}\in \mathcal{J}_{*},$ then we have
\begin{align*}
\langle j_{*}, d(ab) \rangle &=\langle j_{*},(\pi_{*})^{*}( D\circ q(ab)) \rangle \\
&=\langle j_{*},(\pi_{*})^{*} ( D(a+\mathcal{I})(b+ \mathcal{I})) \rangle \\
&=\langle j_{*},(\pi_{*})^{*} (\sigma(a+\mathcal{I})\cdot D(b+\mathcal{I})+D(a+\mathcal{I}).\sigma(b+\mathcal{I}) \rangle \\
&=\langle \pi_{*}( j_{*}),(\sigma(a)+\mathcal{I})\cdot D(b+\mathcal{I})+D(a+\mathcal{I}).(\sigma(b)+\mathcal{I}) \rangle \\
&=\langle \pi_{*}( j_{*})\cdot (\sigma(a)+\mathcal{I}),D(b+\mathcal{I}) \rangle+\langle (\sigma(b)+\mathcal{I})\cdot \pi_{*}( j_{*}), D(a+\mathcal{I}) \rangle \\
&=\langle \pi_{*}( j_{*})\cdot \sigma(a),D(b+\mathcal{I}) \rangle+\langle \sigma(b)\cdot \pi_{*}( j_{*}), D(a+\mathcal{I}) \rangle \\
&=\langle \pi_{*}( j_{*}\cdot \sigma(a), D(b+\mathcal{I}) \rangle+\langle \pi_{*}(\sigma(b)\cdot j_{*}), D(a+\mathcal{I}) \rangle \\
&=\langle j_{*}\cdot \sigma(a),(\pi_{*})^{*}(D(b+\mathcal{I})) \rangle+\langle \sigma(b)\cdot j_{*},(\pi_{*})^{*}( D(a+\mathcal{I})) \rangle \\
&=\langle j_{*},\sigma(a)\cdot (\pi_{*})^{*} (D \circ q(b))+(\pi_{*})^{*} ( D\circ q(a))\cdot \sigma(b) \rangle\\
&=\langle j_{*},\sigma(a)\cdot d(b)+d(a)\cdot \sigma(b) \rangle.
\end{align*}
So there exist an element $\lambda \in \mathcal{J}$ such that $d(a)=\sigma(a)\cdot \lambda - \lambda \cdot \sigma(a) $ ($a \in \mathcal{A}$). Let $m$ be the restriction of $\lambda$ to $\mathcal{I_{*}}$, then $m \in \mathcal{I}$ and for $i_{*}\in \mathcal{I_{*}},$ we have
\begin{align*}
\langle i_{*},\sigma(a)\cdot m -m\cdot \sigma(a) \rangle &=\langle i_{*}\cdot \sigma(a)-a\cdot i_{*},m \rangle\\
&=\langle i_{*}\cdot \sigma(a)-\sigma(a)\cdot i_{*},\lambda \rangle\\
&=\langle i_{*} , \sigma(a)\cdot \lambda -\lambda \cdot \sigma(a) \rangle\\
&=\langle i_{*}, (\pi_{*})^{*}\circ D\circ q(a) \rangle\\
&=\langle \pi_{*}(i_{*}), D\circ q(a) \rangle\\
&=\langle \pi_{*}(i_{*}), D(a+\mathcal{I}) \rangle=0.
\end{align*}
The reason of last equality is that $\pi_{*}$ is the projection on $^{\bot}\mathcal{I}$, so if $i_{*} \in \mathcal{I_{*}}$ since $\mathcal{I_{*}}= \dfrac{\mathcal{A_{*}}}{^{\bot}\mathcal{I}}$ thus $i_{*}$ is not in $^{\bot}\mathcal{I}$ so $\pi_{*}(i_{*})=0$.
Therefore $\sigma(a)\cdot m=m\cdot \sigma(a)$ for each $(a \in \mathcal{A})$. Since $\mathcal{I}$ has the $\sigma-$dual trace extension property, then there exist a $\kappa\in \mathcal{A}$ such that $\kappa \vert_{\mathcal{I_{*}}}=m$ and $ a\cdot \kappa- \kappa \cdot a=0$ ($a \in \mathcal{A}$). Let $\tau$ be the restriction of $\kappa$ to $\mathcal{J_{*}}$. Then $\tau \in \mathcal{J}$ and $\lambda - \tau=0$ on $\mathcal{I_{*}}$. Therefore $\lambda- \tau \in \dfrac{\mathcal{J}}{\mathcal{I}}$.

 \noindent Now let $x$ be in $(\dfrac{\mathcal{J}}{\mathcal{I}})_{*}$, then there exist a $j_{*}\in \mathcal{J}_{*}$ such that $\pi_{*}(j_{*})=x$, so we have
\begin{align*}
\langle x , D(a+\mathcal{I})\rangle &=\langle \pi_{*}(j_{*}), D(a+\mathcal{I})\rangle\\
&=\langle j_{*},\sigma(a)\cdot \lambda - (\sigma(a)\cdot \tau-\tau \cdot a)-\lambda \cdot \sigma(a)\rangle\\
&=\langle j_{*},\sigma(a)\cdot \lambda-\sigma(a)\cdot \tau+\tau \cdot \sigma(a)- \lambda .\sigma(a)\rangle\\
&=\langle j_{*}, \sigma(a)\cdot (\lambda- \tau)-(\lambda-\tau)\cdot \sigma(a). \rangle
\end{align*}
If $j_{*}\in$ $^{\bot}\mathcal{I},$ then by definition of $\pi_{*},$ we have $\pi_{*}(j_{*})= j_{*}.$ Also if $j_{*}$ is not in $^{\bot}\mathcal{I},$ then $\pi_{*}(j_{*})= 0$. In the first case, we have
\begin{align*}
\langle j_{*},\sigma(a)\cdot (\lambda- \tau)-(\lambda-\tau)\cdot \sigma(a) \rangle &=\langle \pi_{*}(j_{*}),\sigma(a)\cdot (\lambda-\tau)-(\lambda-\tau)\cdot \sigma(a)\rangle\\
&=\langle x , \sigma(a)\cdot (\lambda- \tau)-(\lambda-\tau)\cdot \sigma(a) \rangle.
\end{align*}
Hence \begin{center}
$D(a+\mathcal{I})= \sigma(a)\cdot (\lambda- \tau)-(\lambda-\tau)\cdot \sigma(a) .$\end{center}
It means that $D$ is an inner $\sigma$-derivation. In the second case also, we have $D$ is an inner $\sigma$-derivation. So we canclude that $\dfrac{\mathcal{A}}{\mathcal{I}}$ is $\sigma$-ideally Connes amenable.

 \end{proof}



 \begin{proposition}\label{prop1}
\noindent Let $\sigma$ be a $w^{*}-$continuous endomorphism of dual Banach algebra $\mathcal{A}$. If $H^{1}_{\sigma, w^{*}}(\mathcal{A^{\sharp}},\mathcal{A^{\sharp}})=\lbrace 0 \rbrace, $ then $H^{1}_{\sigma, w^{*}}(\mathcal{A},\mathcal{A})=\lbrace 0 \rbrace $.
\end{proposition}
\begin{proof}
\noindent Let $D: \mathcal{A} \longrightarrow \mathcal{A}$ be a $w^{*}-$continuous $\sigma-$derivation. Note that $\mathcal{A}$ is a $\mathcal{A^{\sharp}}-$bimodule with the following module action
\begin{equation}
(a+\alpha )\cdot b = a\cdot b +\alpha b, \quad b\cdot(a+\alpha)= b\cdot a + \alpha b, \nonumber
\end{equation}
for all $a,b \in \mathcal{A}$ and $\alpha \in \mathbb{C}$. Define $\widehat{D}:A^{\sharp} \longrightarrow A^{\sharp} $ with $\widehat{D}(a+\alpha)= D(a).$ Clearly $\widehat{D}$ is a $w^{*}-$continuous $\sigma-$derivation and we can look at it as a function into $\mathcal{A^{\sharp}}$. Since $H^{1}_{\sigma, w^{*}}(\mathcal{A^{\sharp}},\mathcal{A^{\sharp}})=\lbrace 0 \rbrace $, so there exists $t\in \mathcal{A}$ such that $\widehat{D}= \delta_{t}^{\sigma}$. Hence for each $a\in \mathcal{A},$ we have
\begin{eqnarray}
D(a)&=& \widehat{D}(a+\alpha)\nonumber\\
&=& \widehat{\sigma}(a+\alpha)\cdot t - t\cdot \widehat{\sigma}(a+\alpha)\nonumber\\
&=& \sigma(a)\cdot t-t\cdot \sigma(a). \nonumber
\end{eqnarray}
\noindent This shows that $D$ is $\sigma-$inner. Thus $H^{1}_{\sigma, w^{*}}(\mathcal{A},\mathcal{A})=\lbrace 0 \rbrace $.
\end{proof}

 \begin{proposition}\label{prop2}
\noindent Let $\mathcal{A}$ be a dual Banach algebra and $\mathcal{I}$ be a $w^{*}-$closed two-sided ideal of $\mathcal{A}$ with a bounded approximate identity. Suppose that $\sigma$ be an idempotent endomorphism of $\mathcal{A}$ such that $\sigma(\mathcal{I})=\mathcal{I}$. Then $H^{1}_{\sigma, w^{*}}(\mathcal{I},\mathcal{I})=\lbrace 0 \rbrace $ if and only if $H^{1}_{\sigma, w^{*}}(\mathcal{A},\mathcal{I})=\lbrace 0 \rbrace $.
\end{proposition}

 \begin{proof}
\noindent Suppose that $H^{1}_{\sigma, w^{*}}(\mathcal{I},\mathcal{I})=\lbrace 0 \rbrace. $ Let $D:\mathcal{A}\longrightarrow \mathcal{I}$ be $w^{*}-$continuous $\sigma-$derivation and $i: \mathcal{I}\longrightarrow \mathcal{A}$ be the embedding map. Then $d= D\vert_{\mathcal{I}}: \mathcal{I}\longrightarrow \mathcal{I}$ is a $\sigma-$derivation. So there exists $t_{0}\in \mathcal{I}$ such that $d= \delta_{t_{0}}^{\sigma}$. Since $\mathcal{I}$ has a bounded approximate identity and $\sigma(\mathcal{I})=\mathcal{I}$, so $\overline{\sigma(\mathcal{I}^{2})}= \overline{\mathcal{I}^{2}}=\mathcal{I}$. On the other hand, since $\mathcal{I}=\sigma(\mathcal{I})\cdot \mathcal{I}\cdot \sigma(\mathcal{I})$, then we have $\mathcal{I_{*}}=\sigma(\mathcal{I})\cdot \mathcal{I_{*}}\cdot \sigma(\mathcal{I})$. So for $i,j\in \mathcal{I}$ and $i_{*}\in \mathcal{I_{*}},$ we have
\begin{align*}
\langle \sigma(i) i_{*}\sigma(j),& D(a) \rangle = \langle \sigma(i) i_{*}, \sigma(j) D(a) \rangle\\
&=\langle \sigma(i) i_{*}, D(ja)- D(j)\sigma(a)\rangle\\
&=\langle \sigma(i) i_{*}, \sigma(ja)t_{*}-t_{*}\sigma(ja) \rangle - \langle \sigma(i) i_{*}, (\sigma(j)t_{*}- t_{*}\sigma(j))\sigma(a) \rangle\\
&=\langle \sigma(i) i_{*}\sigma(j), \sigma(a)t_{\star} \rangle - \langle \sigma(a) \sigma(i)i_{*}, t_{*}\sigma(j) \rangle\\
&-\langle \sigma(a) \sigma(i)i_{*}, \sigma(j)t_{*} \rangle + \langle \sigma(a) \sigma(i)i_{*}, t_{*}\sigma(j) \rangle\\
&= \langle \sigma(i) i_{*}\sigma(j), \sigma(a) t_{*}-t_{*}\sigma(a) \rangle\\
&=\langle \sigma(i) i_{*}\sigma(j), \delta_{t_{*}}^{\sigma}(a) \rangle.
\end{align*}
Therefore $D=\delta_{t_{*}}^{\sigma}$, and so $D$ is $\sigma-$inner. Conversly, let $H^{1}_{\sigma, w^{*}}(\mathcal{A},\mathcal{I})=\lbrace 0 \rbrace $, and let $D: \mathcal{I}\longrightarrow \mathcal{I}$ be $w^{*}-$continuous $\sigma-$derivation. Since $\mathcal{I}$ is neo-unital Banach $\mathcal{I-}$bimodule i.e. $\mathcal{I}= \sigma(\mathcal{I})\cdot \mathcal{I}\cdot \sigma(\mathcal{I})$, by [\cite{myr}, Proposition 4.14], $D$ has an extension $\widehat{D}: \mathcal{A}\longrightarrow \mathcal{I}$, such that $\widehat{D}$ is also $\sigma-$derivation, now by hypothesis $\widehat{D}$ is $\sigma-$inner. Thus $H^{1}_{\sigma, w^{*}}(\mathcal{I},\mathcal{I})=\lbrace 0 \rbrace $.
\end{proof}

\begin{proposition}\label{lem9}
\noindent Let $\sigma$ be a $w^{*}$-continuous endomorphism of dual Banach algebra $\mathcal{A}$. If $\mathcal{A^{\sharp}} $ is $\widehat{\sigma}-$ideally Connes amenable, then $\mathcal{A}$ is $\sigma-$ideally Connes amenable, where $\widehat{\sigma}$ is defined by $\widehat{\sigma}(a+\alpha)= \sigma(a), \quad a\in \mathcal{A}, \alpha \in \mathbb{C}$.
\end{proposition}

\begin{proof}
\noindent Let $\mathcal{I}$ be $w^{*}$-closed two-sided ideal in $\mathcal{A}$, and $D: \mathcal{A}\longrightarrow \mathcal{I}$ be a
$w^{*}-$continuous $\sigma-$derivation. It is easy to see that $\mathcal{I}$ is a $w^{*}$-closed two-sided ideal of $\mathcal{A^{\sharp}} $, and $\widehat{D}:\mathcal{A^{\sharp}} \longrightarrow \mathcal{I}$ with $\widehat{D}(a+\alpha)= D(a) $, is a $w^{*}-$continuous $\sigma-$derivation. Hence there exists $t\in \mathcal{I}$ such that $\widehat{D}=\delta_{t}^{\sigma}$. So for each $a\in \mathcal{A},$ we have
\begin{eqnarray}
D(a) &=& \widehat{D}(a+\alpha)\nonumber\\
&=&\widehat{\sigma}(a+\alpha)\cdot t-t \cdot \widehat{\sigma}(a+\alpha)\nonumber\\
&=& \sigma(a)\cdot t-t\cdot \sigma(a).\nonumber
\end{eqnarray}

 \noindent This shows that $D$ is $\sigma-$inner. Thus $\mathcal{A}$ is $\sigma-$ideally Connes amenable.
\end{proof}
\begin{proposition}\label{lem10}
\noindent Suppose that $\mathcal{A}, \mathcal{B}$ are dual Banach algebras and $\phi:\mathcal{A}\longrightarrow \mathcal{B}$ be a $w^{*}-$continuous epimorphism. If $\mathcal{A}$ is Connes amenable then $\mathcal{B}$ is $\sigma-$ideal-Connes amenable.
\end{proposition}
\begin{proof}
\noindent Let $\sigma:\mathcal{B}\longrightarrow \mathcal{B}$ be a $w^{*}-$continuous endomorphism on $\mathcal{B}$ and $\mathcal{I}$ be a $w^{*}$-closed two-sided ideal of $\mathcal{B}$, then $\mathcal{I}$ is a normal dual $\mathcal{A-}$bimodule with the following actions given by :\begin{center}
$a\cdot i=\sigma(\phi(a))\cdot i ,\hspace{.5cm}i\cdot a= i\cdot \sigma(\phi(a))\hspace{.5cm}(a\in \mathcal{A}, i\in \mathcal{I}).$
\end{center}
Now let $D: \mathcal{B}\longrightarrow \mathcal{I}$ be a $w^{*}-$continuous $\sigma-$derivation, it is easy to check that $Do\phi: \mathcal{A}\longrightarrow \mathcal{I}$ is a $w^{*}-$continuous derivation, since $\mathcal{I}$ is a normal dual $\mathcal{A-}$bimodule, there exists $t\in \mathcal{I}$ such that $Do\phi(a)= a\cdot t-t\cdot a=\sigma(\phi(a))\cdot t-t\cdot \sigma(\phi(a))$. Since $\phi$ is an epimorphism, so for each $b\in \mathcal{B}$, there exists $a\in \mathcal{A}$ such that $b=\phi(a)$, and we have
\begin{equation}
D(b)=\sigma(b)\cdot t-t\cdot \sigma(b). \nonumber
\end{equation}

 \noindent This shows that $D$ is a $\sigma-$inner derivation. Thus $\mathcal{B}$ is $\sigma-$ideal-Connes amenable.
\end{proof}


\begin{theorem}\label{thm12}
\noindent Suppose that $\mathcal{I}$ is a $w^{*}$-closed two-sided ideal in dual Banach algebra $\mathcal{A}.$ If $\mathcal{I}$ is $\sigma-$Connes amenable and $\dfrac{\mathcal{A}}{\mathcal{I}}$ is Connes-amenable, then $\mathcal{A}$ is $\sigma-$Connes amenable.
\end{theorem}
\begin{proof}
\noindent We follow the method of proof of [\cite{myr}, Proposition 4.18]. Let $X$ be a dual Banach $\mathcal{A}-$bimodule and $D: \mathcal{A}\longrightarrow X$ be a $w^{*}-$continuous $\sigma-$derivation, $X$ is a Banach $\mathcal{I-}$bimodule too. Clearly $d=D\vert_{\mathcal{I}}: \mathcal{I}\longrightarrow X$ is a $\sigma-$derivation and by the $\sigma-$Connes amenability of $\mathcal{I},$ there exists $x_{0}\in X$ such that $D=\delta_{x_{0},w^{*}}^{\sigma}.$ Therefore for each $i\in \mathcal{I},$ we have
\begin{equation}
d(i)=\sigma(i)\cdot x_{0}-x_{0}\cdot \sigma(i). \nonumber
\end{equation}
\noindent Set $D_{1}=D- \delta_{x_{0},w^{*}}^{\sigma}$, clearly $D_{1}$ is $\sigma-$derivation and $D_{1}\vert _{\mathcal{I}}=0$.\\
Now let $X_{0}=\overline{span}^{w^{*}}(X\cdot \sigma(\mathcal{I}))\cup (\sigma(\mathcal{I})\cdot X)$. $\dfrac{X}{X_{0}}$ is a dual Banach $\frac{\mathcal{A}}{\mathcal{I}}-$bimodule via the following actions: \begin{center}
$(a+\mathcal{I})(x+X_{0})=\sigma(a)x+X_{0}, \quad (x+X_{0})(a+\mathcal{I})=x\sigma(a)+X_{0},\hspace{.2cm}(x\in X,a\in \mathcal{A}).$
\end{center}
\noindent Now, we define $\widehat{D}: \dfrac{\mathcal{A}}{\mathcal{I}}\longrightarrow \dfrac{X}{X_{0}}$; $\langle g_{*}, \widehat{D}(a+\mathcal{I}) \rangle= \langle g_{*}, D_{1}(a) \rangle$, for $g_{*}\in (\dfrac{X}{X_{0}})_{*}=^{\perp}X_{0} $.

 \noindent Let $a+\mathcal{I}=a^{\prime}+\mathcal{I}$, for some $a,a^{\prime}\in \mathcal{A}$.
So $a-a^{\prime}\in \mathcal{I}$, and we have $D_{1}(a-a^{\prime})=0$. Thus $D_{1}(a)= D_{1}(a^{\prime})$, so $\widehat{D}(a+\mathcal{I})=\widehat{D}(a^{\prime}+\mathcal{I})$, which shows that $\widehat{D}$ is well defined. We claim that $\widehat{D}$ is a derivation.

 \begin{align*}
\langle g_{*}, \widehat{D}(a+\mathcal{I})(b+\mathcal{I})&=\langle g_{*},D_{1}(ab) \rangle = \langle g_{*},\sigma(a)D_{1}(b)+D_{1}(a)\sigma(b) \rangle \\
&= \langle g_{*}\sigma(a), D_{1}(b) \rangle +\langle \sigma(b) g_{*}, D_{1}(a) \rangle \\
&=\langle g_{*}\cdot (a+\mathcal{I}), \widehat{D}(b+\mathcal{I}) \rangle +\langle (b+\mathcal{I})\cdot g_{*} , \widehat{D}(a+\mathcal{I}) \rangle \\
&=\langle g_{*}, (a+\mathcal{I})\cdot \widehat{D}(b+\mathcal{I}) \rangle +\langle g_{*} , \widehat{D}(a+\mathcal{I})\cdot (b+\mathcal{I}) \rangle.
\end{align*}

 \noindent So there exists $t\in \dfrac{X}{X_{0}}$, such that $\widehat{D}=\delta_{t}$. Thus we have \begin{align*}
\langle g_{*}, D_{1}(a) \rangle &=\langle g_{*} , \widehat{D}(a+\mathcal{I}) \rangle \\
&=\langle g_{*},(a+\mathcal{I})\cdot t-t \cdot (a+\mathcal{I}) \rangle \\
&=\langle g_{*}\cdot (a+\mathcal{I}), t \rangle - \langle (a+\mathcal{I}) \cdot g_{*} , t \rangle \\
&=\langle g_{*}\cdot \sigma(a) , t \rangle - \langle \sigma(a) \cdot g_{*}, t \rangle \\
&= \langle g_{*} , \sigma(a)\cdot t- t\cdot \sigma(a).
\end{align*}

 \noindent So $D_{1}=D-\delta_{x}$, and therefore $D=\delta_{x-t, w^{*}}^{\sigma}.$
\end{proof}
\begin{theorem}\label{pro13}
\noindent Let $\mathcal{A}$ be a $\sigma-$Connes amenable dual Banach algebra. Then $\sigma(\mathcal{A})$ has an identity.
\end{theorem}
\begin{proof}
\noindent The proof is standard, see \cite{run1}. We provide the proof for the benefit of readers. Consider $X=\mathcal{A}$ as a dual Banach $\mathcal{A-}$bimodule with the trivial left action, that is: \begin{center}
$a\cdot x= 0, \hspace{.5cm} x\cdot a= xa, \hspace{.5cm}(a\in \mathcal{A}, x\in X).$
\end{center}
Then $D: \mathcal{A}\longrightarrow X$, defined by $D(a)=\sigma(a)$, is a $w^{*}-$continuous $\sigma-$derivation. Since $\mathcal{A}$ is $\sigma-$Connes amenable, there exists $E\in X$ such that $D=\delta_{E}^{\sigma}$. Hence
\begin{equation}
\sigma(a)= \sigma(a)\cdot E- E\cdot \sigma(a)\quad (a\in \mathcal{A}). \nonumber
\end{equation}

 \noindent Hence $-E$ is an identity element for $\sigma(\mathcal{A})$.
\end{proof}

\begin{proposition}\label{thm19}
\noindent Let $\mathcal{A},\mathcal{B}$ be dual Banach algebras, $\sigma \in Hom^{w^{*}}(\mathcal{B})$ and $\phi: \mathcal{A}\longrightarrow \mathcal{B}$ be a $w^{*}-$continuous epimorphism. If $\mathcal{A}$ is weakly amenable and commutative, then $\mathcal{B}$ is $\sigma-$ideally Connes amenable.
\end{proposition}
\begin{proof}
\noindent Let $\mathcal{I}$ be a $w^{*}-$closed two-sided ideal of $\mathcal{B}$ and $D:\mathcal{B}\longrightarrow \mathcal{I}$ be an arbitrary $w^{*}-$continuous $\sigma-$derivation. Then $\mathcal{I}$ becomes a dual Banach $\mathcal{A-}$bimodule with the following module actions: \begin{center}
$a\cdot i=\sigma(\phi(a))\cdot i, \quad i\cdot a=i \cdot \sigma(\phi(a))\hspace{.2cm}(a\in \mathcal{A},i\in \mathcal{I})$.
\end{center}
The bounded linear mapping $Do\phi: \mathcal{A}\longrightarrow \mathcal{I}$ is a derivation. It is easy to see that $\mathcal{B}$ is commutative and therefore $\mathcal{I}$ is a symmetric Banach $\mathcal{B-}$bimodule. Hence $\mathcal{I}$ is a symmetric Banach $\mathcal{A-}$bimodule. Now $\mathcal{A}$ is weakly amenable, thus $H^{1}(\mathcal{A}, \mathcal{I})=\lbrace 0 \rbrace$. So $Do \phi=0$. Consequently $D=0$, and $\mathcal{B}$ is $\sigma-$ideally Connes amenable.
\end{proof}

 \section{Examples}\label{Se4}
\noindent In this section, we give some examples to illustrate the new notion of $\sigma-$ideally Connes amenability introduced in this work. These examples show that the notion of $\sigma-$ideally Connes amenability is different from ideally Connes amenable. In doing this, we give some examples of $\sigma-$ideally Connes amenable dual Banach algebras that are not ideally Connes amenable.
\example
\noindent Let $\mathcal{A}$ be a dual Banach algebra and $\phi \in$ Ball$\mathcal{A^{*}}$ be a $w^{*}-$continuous homomorphism such that for all $w^{*}-$closed two sided ideal $\mathcal{I}$ of $\mathcal{A}$, $\phi\mid_{\mathcal{I}}\neq 0.$ Then $\mathcal{A}$ with the product
\begin{equation}
a\cdot b=\phi(a)b \quad a,b \in \mathcal{A} \nonumber
\end{equation}
becomes a Banach algebra, since maltiplication is separately $w^{*}-$continuous so $\mathcal{A}$ is a dual Banach algebra. We denote this algebra with $\mathcal{A_{\phi}}$. It is easy to see that $\mathcal{A_{\phi}}$ has a left identity $e$, while it has not right approximate identity.
So by [\cite{mbe}, Proposition 2.3], $\mathcal{A_{\phi}}$ is not ideally Connes amenable. Now suppose that $\sigma: \mathcal{A_{\phi}} \longrightarrow \mathcal{A_{\phi}}$ is defined by
\begin{equation}
\sigma(a)=\phi(a)e. \nonumber
\end{equation}
\noindent We have
\begin{equation}
\sigma^{2}(a) = \sigma(\phi(a)e) \nonumber\\
= \phi(a)\sigma(e) \nonumber\\
= \phi(a)\phi(e)e \nonumber\\
= \phi(a)e \nonumber\\
= \sigma(a). \nonumber
\end{equation}

 \noindent Thus $\sigma$ is idempotent. It is easy to see that $e$ is identity for $\sigma(\mathcal{A_{\phi}})$. Let $\mathcal{I}$ be a $w^{*}-$closed two sided ideal of $\mathcal{A_{\phi}}$ and $d: \mathcal{A_{\phi}}\longrightarrow \mathcal{I}$ be a non-zero $w^{*}-$continuous $\sigma-$derivation. For each $a,b\in \mathcal{A_{\phi}},$
we have
\begin{equation}
d(a\cdot b)=\sigma(a)\cdot d(b)+d(a)\cdot \sigma(b). \nonumber
\end{equation}
\begin{equation}
d(\phi(a)b)=\phi(a)e\cdot d(b)+d(a)\cdot \phi(b)e. \nonumber
\end{equation}
So we have

 \begin{equation}
\phi(a)d(b)=\phi(a)d(b)+\phi(b)d(a)\cdot e. \nonumber
\end{equation}

 \noindent Thus $\phi(b)d(a)\cdot e=0.$ Since $\phi\neq 0,$ we have $d(a)\cdot e=0.$ Thus $\phi(d(a)) e=0$, since $d(a)$ is in $\mathcal{I}$ and $\phi\mid_{\mathcal{I}}\neq 0 $, so we conclude that $e=0$, that is a contradiction. It means that every $\sigma-$derivation is zero, so it is inner. Thus $\mathcal{A_{\phi}}$ is $\sigma-$ideally Connes amenable.

 \example\label{ex2}
\noindent Let $\mathcal A=l^{1} (\mathbb N).$ This algebra was introduced and studied in [\cite{mbe}, Example 2.8]. It is not ideal Connes amenable, since it has no approximate identity. Let $\sigma$ be a $w^{*}-$continuous idempotent endomorphism on $l^{1} (\mathbb N)$ such that for all $a\in l^{1} (\mathbb N)$, $\sigma(a)(1)=a(1)$. If $\mathcal I$ be a weak$^*$-closed two-sided ideal of $\mathcal A$ with $\mathcal I\neq \mathcal A$, we get $\mathcal I\subseteq \{f\in \mathcal A; f(1)=0\}$. For such ideal $\mathcal I$, let $D: \mathcal A\longrightarrow \mathcal I$ be a weak$^*$-continuous $\sigma-$derivation. Set $f \in \mathcal A,$ define the mapping $\widetilde{f}:\mathbb N \longrightarrow \mathbb C$, by $\widetilde{f} (1)=0$ and $\widetilde{f} (n)=f(n)$ for $n\geq 2$. Let $e\in l^{1} (\mathbb N)$ be such that $e(n)=0$ for $n\neq 1$ and $e(1)=1$. Then we find $f=f\cdot e+\widetilde{f}$. We conclude that
\begin{equation}
D(f)=\sigma(f)(1)D(e)+D(\widetilde{f}). \nonumber
\end{equation}

 \noindent For each $g\in \mathcal{A}_*$, we have
\begin{equation}
\langle D(\widetilde{f})\cdot \sigma(e), g\rangle=\langle D(\widetilde{f}), \sigma(e)\cdot g\rangle=\langle D(\widetilde{f}), g\rangle
. \nonumber \end{equation}

 \noindent Since $D(\widetilde{f})\in \mathcal I$, $D(\widetilde{f})(1)=0$, so
\begin{equation}
D(\widetilde{f})\cdot \sigma(e)=D(\widetilde{f})(1) \sigma(e)=0,\nonumber
\end{equation} and we conclude that $D(\widetilde{f})=0.$ Thus,
\begin{equation}
D(f)=\sigma(f)(1)D(e)=\sigma(f)\cdot D(e), \nonumber
\end{equation}
and so
\begin{equation}
D(f)=\sigma(f)\cdot D(e)- D(e) \cdot \sigma(f).\nonumber
\end{equation}

 \noindent In the last equality, since $D(e)\in \mathcal{I}$ and $D(e)(1)=0$, so $D(e) \cdot \sigma(f)=0$. Therefore for any weak$^*$-closed two-sided ideal $\mathcal I$ of $\mathcal A$ with $\mathcal I\neq \mathcal A, H^{1}_{\sigma, w^{*}}(\mathcal{A},\mathcal{I})=\lbrace 0 \rbrace $.

 \noindent We know that $C_{0} (\mathbb N)$ is dense in $l^{1} (\mathbb N)$. Let $a\in l^{1} (\mathbb N)$, then there is a net
$\{a_{\alpha}\}\subset C_{0} (\mathbb N),$ such that $a_{\alpha}\longrightarrow a$ in $w^{*}$-topology. For $f \in C_{0} (\mathbb N)^{*}$, define \begin{equation}
\langle a,f \rangle:=lim-w^{*}\langle a_{\alpha}, f \rangle. \nonumber
\end{equation}
\noindent We can define the left module action of $l^{1} (\mathbb N)$ on $C_{0} (\mathbb N)^{*}$ by
\begin{equation}
a\cdot f= \langle a,f \rangle e. \nonumber
\end{equation}
\noindent We can see that
\begin{equation}
a\cdot (b\cdot f)=(ab)\cdot f, \nonumber
\end{equation}
and
\begin{equation}
\Vert a\cdot f \Vert \leq \Vert a \Vert \Vert f \Vert. \nonumber
\end{equation}
\noindent The right module action is defined by $f\cdot a=a(1)f$. It is easily verified that $C_{0} (\mathbb N)^{*}$ is a $l^{1} (\mathbb N)-$bimodule.

 \noindent Let $D$ be a weak$^*$-continuous $\sigma-$derivation from $l^{1} (\mathbb N)$ to $l^{1} (\mathbb N) \cong C_{0} (\mathbb N)^{*}$. For all $a,b \in l^{1} (\mathbb N),$ we have

 \begin{equation}
D(a.b)=D(a)\cdot \sigma(b)+ \sigma(a)\cdot D(b).\nonumber
\end{equation}

 \noindent Thus
\begin{equation}
a(1)D(b)=\sigma(a)(1)D(a)+\langle \sigma(a), D(b) \rangle e. \nonumber
\end{equation}
\noindent Set $a=b$, then we get
\begin{equation}
\langle \sigma(a), D(a) \rangle=0. \nonumber
\end{equation}

 \noindent So for all $a,b \in l^{1} (\mathbb N),$ we have
\begin{align*}
0=\langle \sigma(a.b), &D(a.b) \rangle=\langle \sigma(a.b) , D(a)\cdot \sigma(b)+ \sigma(a)\cdot D(b) \rangle\\
&=\langle \sigma(a).\sigma(b), \sigma(b)(1)D(a)+ \langle \sigma(a), D(b) \rangle e \rangle\\
&= \langle \sigma(a).\sigma(b), \sigma(b)(1)D(a)+ \langle \sigma(a), D(b) \rangle \langle \sigma(a)\sigma(b), e \rangle\\
&= \sigma(b)(1) \langle \sigma(a).\sigma(b), D(a) \rangle + \sigma(a)(1)\sigma(b)(1)\langle \sigma(a), D(b) \rangle\\
&= \sigma(b)(1).\sigma(a)(1)\langle \sigma(b), D(a) \rangle + \sigma(a)(1)\sigma(b)(1)\langle \sigma(a), D(b) \rangle.
\end{align*}
So we conclude that
\begin{equation}
\langle \sigma(a), D(b) \rangle = -\langle \sigma(b), D(a) \rangle. \nonumber
\end{equation}

 \noindent Let $t\in \sigma(\mathcal{A})$, then there is a, $b\in l^{1} (\mathbb N) $ such that $t=\sigma(b)=\sigma(\sigma(b)).$ We obtain
\begin{align*}
\langle t, D(a) \rangle&=\langle t, D(e.a) \rangle\\
&=\langle t, D(e)\cdot \sigma(a)+ \sigma(e)\cdot D(a) \rangle\\
&= \langle t, D(e)\cdot \sigma(a)\rangle +\langle t, \sigma(e)\cdot D(a) \rangle\\
&=\langle \sigma(a).t, D(e) \rangle + \langle t.\sigma(e), D(a) \rangle\\
&=\langle \sigma(a). \sigma(\sigma(b)), D(e) \rangle + \langle \sigma(\sigma(b)).\sigma(e), D(a) \rangle\\
&=\langle \sigma(a). \sigma(b), D(e) \rangle +\langle \sigma(\sigma(b).e), D(a) \rangle\\
&=\langle \sigma(a). \sigma(b), D(e) \rangle - \langle \sigma(a), D(\sigma(b).e) \rangle\\
&= \langle \sigma(a). \sigma(b), D(e) \rangle - \langle \sigma(a), \sigma(b)(1)D(e) \rangle\\
&=\langle \sigma(b), D(e) \cdot \sigma(a)\rangle -\langle \sigma(a), D(e)\cdot \sigma(b) \rangle\\
&= \langle \sigma(b), D(e) \cdot \sigma(a)\rangle -\langle \sigma(b).\sigma(a), D(e) \rangle\\
&= \langle \sigma(b), D(e) \cdot \sigma(a)\rangle -\langle \sigma(b), \sigma(a)\cdot D(e) \rangle\\
&=\langle t, D(e) \cdot \sigma(a)- \sigma(a)\cdot D(e) \rangle.
\end{align*}
So,
\begin{equation}
D(a)=D(e) \cdot \sigma(a)- \sigma(a)\cdot D(e)=\delta_{-D(e)}^{\sigma}(a). \nonumber
\end{equation}
\noindent Thus, we conclude that $ l^{1} (\mathbb N)$ is $\sigma-$ideally Connes amenable.

\example\label{ex3}

 \noindent Let $\mathcal{A}$ be a non-ideally Connes amenable Banach algebra with a right [left] approximate identity. It is known that $\mathcal{A}^{\sharp}$ is not ideally Connes amenable. Consider the map $\sigma:\mathcal{A}^{\sharp} \longrightarrow \mathcal{A}^{\sharp} $ defined by
\begin{equation}
\sigma(a+\lambda e_{\mathcal{A}^{\sharp}})=\lambda, \quad (a\in \mathcal{A}, \lambda \in \mathbb{C}). \nonumber
\end{equation}
Then it is routinely checked that $\sigma$ is weak$^*$-continuous. Let $\mathcal{I}$ be a $w^{*}$-closed two-sided ideal in $\mathcal{A}^{\sharp}$, we show that every $\sigma-$derivation $D:\mathcal{A}^{\sharp} \longrightarrow \mathcal{I} $ is zero. Let $(e_{i})_{i}$ be a right [left] approximate identity for $\mathcal{A}$. A simple calculation shows that $D(ae_{i})=0 $, for each $a\in \mathcal{A}$ and $i$, and consequently $D(a)=0$. Hence
\begin{equation}
D(a+\lambda e_{\mathcal{A}^{\sharp}})=D(a)+\lambda D(e_{\mathcal{A}^{\sharp}})=0. \nonumber
\end{equation}
That is $D=0$ and so $\mathcal{A}^{\sharp}$ is $\sigma$-ideally Connes amenable.

 \section{Conclusion}\label{Se5}
\noindent In this work, we introduced and studied the notions of $\sigma-\mathcal{I}-$ Connes amenability and $\sigma-$ ideally Connes amenability for a dual Banach algebra $\mathcal{A},$ where $\mathcal{I}$ is a $weak^{*}$-closed
two-sided ideal in $\mathcal{A}$ and $\sigma$ is a $weak^{*}$-continuous endomorphism on $\mathcal{A}.$ Under some conditions, we proved
that the ideal Connes-amenability of a dual Banach algebra $\mathcal{A}$ is equivalent to its $\sigma-$ ideally Connes amenability.
Some general theory and hereditary properties on the notion of $\sigma-$ ideally Connes amenability for dual Banach algebras are established. Finally, we presented some useful and important examples to illustrate our results. The results obtained in this work complement and extend some existing results in the literature.\\\\

\bigskip


\end{document}